\documentclass[11pt]{article}
\usepackage{amsmath}
\usepackage{amsthm,amscd,amsfonts}
\usepackage{amssymb, upref, color}
\usepackage{amsmath,amssymb,amsthm}
\usepackage{color}
\usepackage[colorlinks]{hyperref}
\usepackage{graphicx}
\usepackage{epsf,epsfig,subfigure, verbatim}
\usepackage{latexsym,bm}
\usepackage{enumerate}
\usepackage{listings}
\usepackage{cite}
\usepackage{mathrsfs}

\usepackage{color}

\newtheorem{exam}{Example}[section]
\newtheorem{theorem}[exam]{Theorem}
\newtheorem{lemma}[exam]{Lemma}
\newtheorem{remark}[exam]{Remark}

\newtheorem{definition}[exam]{Definition}

\begin{document}
\date{}


\title{ Exponential trichotomy and global linearization of non-autonomous differential equations
}
\author
{
Chaofan Pan$^{a}$\,\,\,\,\,
Manuel Pinto$^{b}$\,\, \,\,
Yonghui Xia$^{a}\footnote{Corresponding author. Yonghui Xia, xiadoc@163.com;yhxia@zjnu.cn. }$
\\
{\small \textit{$^a$ College of Mathematics and Computer Science,  Zhejiang Normal University, 321004, Jinhua, China}}\\
{\small \textit{$^b$ Departamento de Matem\'aticas, Universidad de Chile, Santiago, Chile }}\\
{\small Email: pancf0530@zjnu.edu.cn;  pintoj.uchile@gmail.com; xiadoc@163.com; yhxia@zjnu.cn.}
}

\maketitle

\begin{abstract}

Hartman-Grobman theorem was initially extended to the non-autonomous cases by Palmer.
Usually, dichotomy is an essential condition of Palmer's linearization theorem. Is Palmer's linearization theorem valid for the systems with trichotomy?
In this paper, we obtain  new versions of the linearization theorem if linear system admits  exponential trichotomy on $\mathbb{R}$. {
Furthermore, the equivalent function $\mathscr H(t,x)$ and its inverse $\mathscr L(t,y)$ of our linearization theorems are  H\"{o}lder continuous}. In addition, if a system is periodic, we find the equivalent function $\mathscr H(t,x)$ and its inverse $\mathscr L(t,y)$ of our linearization theorems do not have periodicity or asymptotical periodicity.  To the best of our knowledge, this is the first paper studying the linearization with exponential trichotomy.

\noindent{\bf Keywords:} Exponential trichotomy; Linearization; Periodic;\\

\noindent {\bf MSC 2022:} 34D09;34D10

\end{abstract}
\section{Introduction}
\subsection{Brief history on trichotomy}
In 1930, Perron \cite{Perron} proposed the notion of  (uniform or classical) exponential dichotomy. Later, many dichotomies were introduced, such as nonuniform exponential dichotomy (see Barreira and Valls \cite{Valls1,Valls2}), $(h,k)$-dichotomy (see Naulin and Pinto \cite{Naulin},  Fenner and Pinto \cite{J.Fenner3}), algebraic dichotomy (see Lin \cite{Lin}), $(h,k,\mu,\nu)$-dichotomy (see Zhang et al. \cite{J.Zhang}) and so on. In 1975, Sacker and Sell \cite{SS1} proposed the concept of trichotomy for linear differential systems, decomposing $\mathbb{R}^n$ into stable, unstable and neutral subspaces. Later,  Elaydi and Hajek \cite{EH2} introduced  a stronger notion of trichotomy. Hong, Obaya and Gilet \cite{Hong} considered the existence of a class of ergodic solutions for some differential equations by using exponential trichotomy. Barreira and Valls \cite{Barreira1,Barreira2}  showed that the existence of a nonuniform exponential trichotomy under sufficiently small $C^{1}$ perturbations. Popa, Ceausu and Bagdaser \cite{Popa}  considered   linear discrete-time systems by generalized exponential trichotomy. Adina and  Bogdan \cite{Sasu} study the uniform exponential trichotomy of variational difference equations.  In Banach spaces, Kovacs \cite{Kovacs} considered three concepts of uniform exponential trichotomy on the half-line in the general framework of evolution operators.

\subsection{$C^0$ linearization of the differential equations}
On the other hand,  we are interested in the linearization of the ordinary differential equations. Hartman and Grobman \cite{Hartman, Grobman} made a basic contribution to the linearization problem for autonomous differential equations (called Hartman-Grobman theorem). Later, Hartman-Grobman theorem are generalized in scalar reaction-diffusion equations, Cahn-Hilliard equation, phase field equations and random dynamical systems (see Lu \cite{Lu1}, Bates and Lu \cite{Bates1}, Barreira and Valls \cite{L.Barreira8}). Pugh \cite{C.Pugh} used certain powerful functional analytic skills  to obtain another proof way of Hartman-Grobman theorem. In Banach spaces, Hein and Pr\"{u}ss \cite{Hein-Pruss1} extended Hartman-Grobman theorem to abstract semilinear hyperbolic evolution equations. Reinfelds \cite{A.Reinfelds1} proved that some specific differential equations  are strictly dynamically equivalent. Reinfelds and Sermone  \cite{Reinfelds2} gave a linearization result in nonlinear differential equations with impulse effect.  For the dynamical equivalence of quasilinear impulsive equations, one can refer to Reinfelds \cite{Rein1997,Rein2000},  Sermone \cite{Sermone3,Sermone4},  Reinfelds and $\breve{S}$teinberga \cite{Reinfelds-IJPAM,Reinfelds3}.

In 1973, Palmer \cite{Palmer2} successfully generalized the Hartman-Grobman theorem to non-autonomous differential equations
\begin{equation}
y'=A(t)y+f(t,y). \label{equ1.1}
\end{equation}
In order to weaken the conditions of Palmer's linearization theorem, some improvements were reported: without  exponential dichotomy (see Backes, Dragi\v{c}evi\'{c} and Palmer \cite{BDK-JDE}), for nonuniform dichotomy (see Barreira and Valls \cite{L.Barreira2, L.Barreira3, L.Barreira4, L.Barreira9}), for generalized dichotomy and ordinary dichotomy (Jiang \cite{L.Jiang1,L.Jiang2}), 
for nonuniform contraction (see Casta\~{n}eda and Huerta \cite{Huerta2,Huerta1}), for differential equations with piecewise constant argument (see Zou, Xia and Pinto \cite{Zou-Xia}), for dynamic systems on time scales (Xia et al. \cite{Xia-JDE}, P\"otzche \cite{Potzche1}), for the instantaneous impulsive system (see Fenner and Pinto \cite{J.L.Fenner2}, Xia and Chen \cite{Xia1}), Papaschinopoulos \cite{Papaschinopoulos-A},  Casta\~{n}eda, Gonz\'{a}lze and Robledo \cite{Pablo}, Pinto and Robledo \cite{Pinto-A}, for nonuniform $(h,k,\mu,\nu)$-dichotomy with ordinary differential equations (see Zhang, Fan and Zhu \cite{J.Zhang}), for nonuniform $(h,k,\mu,\nu)$-dichotomy with nonautonomous impulsive differential equations (see Zhang, Chang and Wang \cite {Zhang2}), for non-instantaneous impulsive nonautonomous (see Li, Wang and O'Regan \cite{Li-Wang1,Li-Wang2,Wang1,Wang-IMA}), for admissibility and roughness of nonuniform exponential dichotomies (see Zhou and Zhang \cite{ZhouLF1,ZhouLF2}),  Dragi\v{c}evi\'{c}, Zhang and Zhou \cite{BDX} (admissibility and nonuniform exponential dichotomies), for generalized exponential dichotomies with invariant manifolds (see Zhang \cite{Zhang3}). Above mentioned works are for the $C^0$ linearization. Recently, there are some interesting advance in the $C^1$ linearization for hyperbolic diffeomorphisms (see e.g Backes and  Dragi\v{c}evi\'{c} \cite{BD2}; Dragi\v{c}evi\'{c}, Zhang and Zhang \cite{DZZ-MZ,DZZ-PLMS}; Zhang, Zhang and Jarczyk \cite{ZWM2}; Zhang and Zhang \cite{ZWM3,ZWM4};
Zhang, Lu and  Zhang \cite{ZWM5}).

\subsection{Motivation and novelty}
 Palmer's linearization theorem \cite{Palmer2} needs two essential conditions: (i) the nonlinear term $f$ is bounded and Lipschitzian; (ii) the linear system
\begin{equation}
 x'(t)=A(t)x(t) \label{eq1}
\end{equation}admits  exponential dichotomy.
In this paper,  we  pay particular attention to the effect of the exponential trichotomy imposing on the {\bf linearization of the non-autonomous} ordinary differential equations. Motivated by the works of Palmer \cite{Palmer2}, Backes, Dragi\v{c}evi\'{c} and Palmer \cite{BDK-JDE}, Elaydi and Hajek's exponential  trichotomy (see Elaydi and Hajek  \cite{EH2}), we give  new versions of the linearization theorems based on  exponential trichotomy. The main contributions of the present paper is to improve Palmer's linearization theorem in four aspects:

\noindent {\bf(I)}: The linear system admits  exponential trichotomy, which is weaker than  exponential dichotomy.

\noindent {\bf(II)}: We prove that the equivalent functions $\mathscr H(t,x)$ and its inverse $\mathscr L(t,y)$ are H\"{o}lder continuous.

\noindent {\bf(III)}: The periodicity of the equivalent function is investigated. We prove that periodicity is not an invariant property under exponential trichotomy. More specially, for the periodic systems,  if the linear system admits exponential trichotomy, the equivalent function $\mathscr H(t,x)$ and its inverse $\mathscr L(t,y)$ are not periodic, asymptotically periodic or almost periodic. While, if the linear system admits exponential dichotomy, its equivalent functions are periodic (see \cite{Xia-JDE}).

\noindent {\bf(IV)}:  The nonlinear term $f$ could be unbounded or non-Lipschitzian in our second linearization theorem. 


\subsection{Outline of this paper}
The structure of our paper  as follows.  In Section 2, we give  some basic definitions. In Section 3, we give the our theorems. In Section 4, we  prove our results. Finally, we give some examples to show our linearization theorems.

\section{Statement of main results}
Let $(X,|\cdot|)$ be an arbitrary Banach space. $A(t)$ is a $n\times n$ continuous and bounded matrice
defined on $\mathbb{R}$ respectively. $f:\mathbb{R}\times X\to X$ is a continuous map. \\
\indent Consider the systems
\begin{equation}\label{sys2.1}
x'= A(t)x+f(t,x),
\end{equation}
and
\begin{equation}\label{sys2.2}
x'= A(t)x.
\end{equation}
\begin{definition}\cite{Xia2}\label{def2.1}
Suppose that there exists a function $\mathscr H : \mathbb{R} \times X \to X $ such that\\
$(i)$ for each fixed $t$, $\mathscr H(t,\cdot)$ is a homeomorphism of $X$ into $X$;\\
$(ii)$ $|\mathscr H(t, x)-x|$ is uniformly bounded with respect to $t$;\\
$(iii)$ assume that $\mathscr L(t, \cdot) =\mathscr H^{-1}(t, \cdot)$ also has property $(ii)$;\\
$(iv)$ if $x(t)$ is a solution of system (\ref{sys2.1}), then $\mathscr H(t, x(t))$ is a solution of system (\ref{sys2.2}); and if $y(t)$ is a solution of system (\ref{sys2.2}), then $\mathscr L(t, y(t))$ is a solution
of system (\ref{sys2.1}).\\
\indent If such a map $\mathscr H_{t}(\mathscr H_{t} :=\mathscr H(t, \cdot))$ exists, then  system (\ref{sys2.1}) is topologically conjugated to  system (\ref{sys2.2}) and the transformation $\mathscr H(t, x)$ is called an equivalent function.
\end{definition}
\begin{definition}\cite{EH2}\label{def2.2}
Linear system (\ref{sys2.2}) is said to possess an exponential trichotomy, if there exists  projections $P$, $Q$ and constants $\beta \geq1 $, $\alpha > 0$ such that
\begin{equation}\label{for2.3}
\begin{split}
\begin{cases}
PQ=QP,\quad P+Q-PQ=I,\\
|U(t)PU^{-1}(s)|\leq \beta e^{-\alpha(t-s)} \quad(0\leq s\leq t),\\
|U(t)(I-P)U^{-1}(s)|\leq \beta e^{-\alpha(s-t)}  \quad(t\leq s, s\geq0),\\
|U(t)QU^{-1}(s)|\leq \beta e^{-\alpha(s-t)} \quad(t\leq s\leq 0),\\
|U(t)(I-Q)U^{-1}(s)|\leq \beta e^{-\alpha(t-s)}  \quad(s\leq t, s\leq0),
\end{cases}
\end{split}
\end{equation}
hold; here $U(t)$ is  a  fundamental matrix of the linear system (\ref{sys2.2}).
\end{definition}

\begin{remark}
If we take $P=I-Q$ in the  Definition \ref{def2.2}, then (\ref{for2.3}) becomes
\begin{equation}\label{for2.5}
\begin{split}
\begin{cases}
|U(t)PU^{-1}(s)|\leq \beta e^{-\alpha(t-s)} \quad( s\leq t),\\
|U(t)QU^{-1}(s)|\leq \beta e^{-\alpha(s-t)} \quad(t\leq s).\\
\end{cases}
\end{split}
\end{equation}
We obtain an exponential dichotomy on $\mathbb{R}$.
\end{remark}
\begin{remark}
The first inequality of (\ref{for2.3}) can be divided into the first and the fourth inequalities of (\ref{for2.3}). The second inequality of (\ref{for2.3}) can be divided into the second and the third inequalities of (\ref{for2.3}). Thus,
it is always true that exponential dichotomy on $\mathbb{R}$ implies exponential trichotomy. However, the converse is clearly false as it may be shown
by simple example.
\end{remark}
Next example shows that the linear system admits an exponential trichotomy, but it does not admit an exponential dichotomy.
\begin{exam}\cite{EH2}\label{exam2.5}
\rm{Consider the scalar equation}
\begin{equation}\label{sys2.6}
x'=(-\frac{e^{t}-e^{-t}}{e^{t}+e^{-t}})x.
\end{equation}
Then $x(t)=\frac{2}{e^{t}+e^{-t}}x_{0}$ is the solution of equation (\ref{sys2.6}) with $x(0)=x_{0}$. Now we
take $P=I$, $Q=I$. Obviously, $I-P=0$, $I-Q=0$ , $PQ=QP$, $P+Q-QP=I$.
\begin{align*}
&|x(t)Px^{-1}(s)|=\frac{e^{s}+e^{-s}}{e^{t}+e^{-t}}
\leq\frac{e^{t+s}+e^{t-s}}{e^{2t}}
\leq e^{s-t}+e^{-s-t}
\leq 2e^{-(t-s)}, (t\geq s \geq0),\\
&|x(t)Qx^{-1}(s)|=\frac{e^{s}+e^{-s}}{e^{t}+e^{-t}}
\leq\frac{e^{s-t}+e^{-s-t}}{e^{-2t}}
\leq e^{s+t}+e^{-(s-t)}
\leq 2e^{-(s-t)},(t\leq s \leq0).
\end{align*}
The last two inequalities of definition \ref{def2.2} obviously hold in this case. This implies that equation (\ref{sys2.6}) admits an exponential trichotomy with
\begin{align*}
\beta=2,\quad \alpha=1.
\end{align*}
However, equation (\ref{sys2.6}) doesn't satisfy exponential dichotomy.
\end{exam}
Exponential trichotomy in Definition \ref{def2.2} has the Green function:
\begin{equation}\label{for2.4}
\begin{split}
G(t,s)=
\begin{cases}
U(t)PU^{-1}(s), \quad(0\leq s\leq t),\\
-U(t)(I-P)U^{-1}(s),  \quad(t\leq s, s\geq0),\\
-U(t)QU^{-1}(s), \quad(t\leq s\leq 0),\\
U(t)(I-Q)U^{-1}(s),\quad(s\leq t, s\leq0).
\end{cases}
\end{split}
\end{equation}

Now we consider exponential trichotomy by decomposition of fundamental matrix $U(t)$.
We assume that $U(t)=(\varphi_{1}(t),\varphi_{2}(t),...,\varphi_{n}(t))$ is a  fundamental matrix of the linear system $x'= A(t)x$. $\varphi_{i}(t)$ is bounded on $\mathbb{R}$ for $i=1,...,r-k,r+l,...,n$, $\varphi_{j}(t)$ is unbounded on $t\leq k$ (and bounded on $t\geq k$) for $j=r-k+1,...,r+l-1$.
Then, we chose projection $P=P_{1}+P_{2}+P_{3}=diag\{I_{r-k},0\}+diag\{0_{r-k},I_{r},0\}+diag\{0_{r},I_{r+l-1},0\}$ and  $Q=P_{2}+P_{3}+P_{4}=diag\{0_{r-k},I_{r},0\}+diag\{0_{r},I_{r+l-1},0\}+diag\{0_{r-l+1}, I_{n}\}.$ Furthermore,
\begin{align*}
&PQ=P_{2}+P_{3}=QP,\\
&P+Q-PQ=I.
\end{align*}
We verify projections $P,Q$ satisfy the first inequality of Definition \ref{def2.2}. Then $|U(t)P_{1}U ^{-1}(s)|$  and $|U(t)P_{4}U^{-1}(s)|$ are bounded on $\mathbb{R}$. $|U(t)(P_{2}+P_{3})U^{-1}(s)|$ is unbounded on $t\leq s$ (and bounded on $t\geq s$). We denote that
\begin{align}
|U(t)P_{1}U ^{-1}(s)|&\leq \kappa_{1},\label{fora3.7}\\
|U(t)P_{4}U^{-1}(s)|&\leq \kappa_{2}.\label{fora3.8}
\end{align}
Then we  obtain another  Green function $\tilde{G}(t,s)$,
\begin{equation}\label{foru2.7}
\begin{split}
\tilde{G}(t,s)=
\begin{cases}
U(t)P_{1}U^{-1}(s), \quad(s \leq t),\\
U(t)(P_{2}+P_{3})U^{-1}(s),  \quad(0\leq s\leq t),\\
-U(t)P_{4}U^{-1}(s), \quad(t\leq s),
\end{cases}
\end{split}
\end{equation}
which appears in the statement of Theorem \ref{thmh}. We will study linearization based on these two Green functions.
\section{Main results}
\begin{theorem}\label{thm2.6}
Suppose that  linear system (\ref{sys2.2}) admits an exponential trichotomy (\ref{for2.4}) (that is, fundamental matrix $U(t)$ satisfying (\ref{for2.3}))  and $f(t, x)$ satisfies
\begin{equation}\label{for2.7}
\begin{split}
\begin{cases}
|f(t,x)|\leq \mu,\\
|f(t,x_{1})-f(t,x_{2})|\leq{\gamma}|x_{1}-x_{2}|,\\
3\beta\gamma\alpha^{-1}<1.
\end{cases}
\end{split}
\end{equation}
Then, nonlinear system (\ref{sys2.1}) is topologically conjugated to its linear system (\ref{sys2.2}).
\end{theorem}
\begin{remark}
In Theorem \ref{thm2.6}, if exponential trichotomy reduce to exponential dichotomy, Theorem \ref{thm2.6} still holds. Indeed, that is Palmer's linearization theorem, (see Palmer \cite{Palmer2}).
\end{remark}

In what follows, we introduce an assumption motivated by Backes, Dragi\v{c}evi\'{c} and Palmer \cite {BDK-JDE}.\\
\noindent{\bf Condition (I)}: suppose that there is a continuous function $\Delta_{1}(t,s)> 0$ such that if $z_{1}(t)$ and $z_{2}(t)$ are the solution of system (\ref{sys2.2}), then $|z_{1}(t)-z_{2}(t)|\leq \Delta_{1}(t,s)|z_{1}(s)-z_{2}(s)|$;
  there is another continuous function $\Delta_{2}(t,s)> 0$ such that if $z_{3}(t)$ and $z_{4}(t)$ are the solution of system (\ref{sys2.1}), then $|z_{3}(t)-z_{4}(t)|\leq \Delta_{2}(t,s)|z_{3}(s)-z_{4}(s)|$ ($\Delta_{1}(t,s)>0$ and $\Delta_{2}(t,s)>0$ are continuous functions).
\begin{remark}
Condition (I) is valid, one can refer to Appendix A in \cite {BDK-JDE} for the detail of functions $\Delta_{1}(t,s)$ and $\Delta_{2}(t,s)$.
\end{remark}

\begin{theorem}\label{thm2.7}
Suppose that the conditions in Theorem \ref{thm2.6} and condition ($I$) are satisfied. Let $p>0$ and $0<q<1$ such that
\begin{align*}
m\sup_{t\in\mathbb{R}}\int_{-\infty}^{+\infty}|G(t,s)|\Delta^{q}_{1}(t,s)ds\leq\frac{p}{1+p},\\
m\sup_{t\in\mathbb{R}}\int_{-\infty}^{+\infty}|G(t,s)|\Delta^{q}_{2}(t,s)ds\leq\frac{p}{1+p},
\end{align*}
where $m=\max\{\gamma,2\mu^{1-q}\gamma^q\}$. Then
\begin{equation}
\begin{split}
|\mathscr H(t,x)-\mathscr H(t,x')|&\leq (1+p)|x-x'|^{q}, \quad \mathrm{if} \ |x-x'|<1,\\
|\mathscr L(t,y)-\mathscr L(t,y')|&\leq (1+p)|y-y'|^{q}, \quad  \mathrm{if} \ |y-y'|<1,
\end{split}
\end{equation}
where $\mathscr L(t, \cdot) = \mathscr H^{-1}(t, \cdot)$ and $p$ is a positive constant.
\end{theorem}


\begin{theorem}\label{thm2.8}
Suppose that  linear system (\ref{sys2.2}) admits an exponential trichotomy (\ref{foru2.7}) (that is, fundamental matrix $U(t)$ satisfying (\ref{for2.3}))  and $f(t, x)$ satisfies
\begin{equation}\label{fora2.8}
\begin{split}
\begin{cases}
|f(t,x)|\leq \phi(t),\\
|f(t,x_{1})-f(t,x_{2})|\leq\psi(t)|x_{1}-x_{2}|,\\
\displaystyle\int_{-\infty}^{+\infty}\phi(t)dt<c_{1},\\
\displaystyle\int_{-\infty}^{+\infty}\psi(t)dt<c_{2},
\end{cases}
\end{split}
\end{equation}
where $\phi,\psi:\mathbb{R}\to[0,\infty)$ are integrable functions; $c_{1},c_{2}$ are positive constant. Then, nonlinear system (\ref{sys2.1}) is topologically conjugated to its linear system (\ref{sys2.2}).
\end{theorem}

\begin{theorem}\label{thmh}
Suppose that the conditions in Theorem \ref{thm2.8} and condition (I) are satisfied.
\begin{align*}
\sup_{t\in\mathbb{R}}\int_{-\infty}^{+\infty}|\tilde{G}(t,s)|M(s)\Delta_{1}^{q}(t,s)ds\leq\frac{p}{1+p},\\ \sup_{t\in\mathbb{R}}\int_{-\infty}^{+\infty}|\tilde{G}(t,s)|M(s)\Delta_{2}^{q}(t,s)ds\leq\frac{p}{1+p},
\end{align*}
where $M(t)=\max\{\psi(t),2\phi^{1-q}(t)\psi^{q}(t)\}$ and $q$ is in (0,1). Then
\begin{equation}
\begin{split}
|\mathscr H(t,x)-\mathscr H(t,x')|\leq (1+p)|x-x'|^{q}, \quad \mathrm{if} \ |x-x'|<1,\\
|\mathscr L(t,y)-\mathscr L(t,y')|\leq (1+p)|y-y'|^{q'}, \quad  \mathrm{if} \ |y-y'|<1,
\end{split}
\end{equation}
where $\mathscr L(t, \cdot) = \mathscr H^{-1}(t, \cdot)$ and $p$ is a positive constant.
\end{theorem}


\begin{theorem}\label{thm2.10}
In system (\ref{sys2.1}), assuming that $A(t)$ and $f(t,x)$  has $T$ period with respect to $t$.  Then the equivalent function $\mathscr H(t,x)$ and its inverse $\mathscr L(t,y)$ in  Theorem \ref{thm2.6} and Theorem \ref{thm2.8} do not have  periodicity, asymptotical periodicity or almost periodicity with respect to t.
\end{theorem}

\section{Proofs of main results}
Let $x(t,t_{0},x_{0})$ be the solution of system (\ref{sys2.1}) satisfies the initial value condition $x(t_{0})=x_{0}$, $y(t,t_{0},y_{0})$ be the solution of system (\ref{sys2.2}) satisfies the initial value condition $y(t_{0})=y_{0}$.\\
{\bf Proof  of Theorem \ref{thm2.6}.} \\
{\bf Step 1.} We prove the case of $t\geq0$. Let $\Omega$ denote the space of all continuous maps $h$ satisfies
\begin{align*}
\Omega:=\{h(t)|h:\mathbb{R}\to X, |h(t)| \leq 3\beta\mu\alpha^{-1}\}.
\end{align*}
Then, $(\Omega,\|\cdot\|)$ is a Banach space. For $t\geq0$, $\xi,\eta\in X$ and $h_{1}\in\Omega$, we define the following mapping:
\begin{equation}\label{for3.1}
\begin{split}
&{\mathscr T}h_{1}(t,\eta)=\tilde{h}_{1}(t,\eta)\\
=&\int_{-\infty}^{0}U(t)(I-Q)U^{-1}(s)f(s,y(s,t,\eta)+h_{1}(s,y(s,t,\eta)))ds\\
&+\int_{0}^{t}U(t)PU^{-1}(s)f(s,y(s,t,\eta)+h_{1}(s,y(s,t,\eta)))ds\\
&-\int_{t}^{+\infty}U(t)(I-P)U^{-1}(s)f(s,y(s,t,\eta)+h_{1}(s,y(s,t,\eta)))ds\\
=&\int_{-\infty}^{+\infty}G(t,s)f(s,y(s,t,\eta)+h_{1}(s,y(s,t,\eta))ds,
\end{split}
\end{equation}
where $G(t,s)$ is defined in (\ref{for2.4}). It follows from (\ref{for2.3}), we obtain
\begin{equation}\label{for3.2}
\begin{split}
&|\tilde{h}_{1}(t,\eta)|\\
\leq& \int_{-\infty}^{0}\beta e^{-\alpha(t-s)}\mu ds
+\int_{0}^{t}\beta e^{-\alpha(t-s)}\mu ds + \int_{t}^{+\infty}\beta e^{-\alpha(s-t)}\mu ds\\
\leq&\beta\mu\alpha^{-1}+\beta\mu\alpha^{-1}+\beta\mu\alpha^{-1}\\
\leq&3\beta\mu\alpha^{-1}.
\end{split}
\end{equation}
Moreover,  by differentiating (\ref{for3.1}), we get
\begin{equation}\label{for3.3}
\begin{split}
\tilde{h}^{'}_{1}(t,\eta)=A(t)\tilde{h}_{1}(t,\eta)+f(t,\eta+h_{1}(t,\eta)).
\end{split}
\end{equation}
Hence, from (\ref{for3.2}) and (\ref{for3.3}), we know $\tilde{h}_{1}(t,\eta)$ is continuous and $\tilde{h}_{1}(t,\eta)\in \Omega$.  For any $h_{2},h_{3}\in\Omega$,
from (\ref{for2.7}) and (\ref{for3.1}), we have
\begin{equation}\label{for3.4}
\begin{split}
&|\tilde{h}_{2}(t)-\tilde{h}_{3}(t)|\\
\leq& \int_{-\infty}^{0}\beta e^{-\alpha(t-s)}\gamma|h_{2}(s)-h_{3}(s)| ds+\int_{0}^{t}\beta e^{-\alpha(t-s)}\gamma|h_{2}(s)-h_{3}(s)| ds\\
&+\int_{t}^{+\infty}\beta e^{-\alpha(t-s)}\gamma|h_{2}(s)-h_{3}(s)| ds\\
\leq&3\beta\gamma\alpha^{-1}\|h_{2}-h_{3}\|.
\end{split}
\end{equation}
Note that $3\beta\gamma\alpha^{-1}<1$. Thus, $\mathscr T$ : $\Omega\to\Omega$ is a contraction map. Therefore, there exists a  unique fixed point $h_{1}\in\Omega$ such that
\begin{equation}\label{for3.5}
\begin{split}
&h_{1}(t,\eta)\\
=&\int_{-\infty}^{0}U(t)(I-Q)U^{-1}(s)f(s,y(s,t,\eta)+h_{1}(s,y(s,t,\eta)))ds\\
&+\int_{0}^{t}U(t)PU^{-1}(s)f(s,y(s,t,\eta)+h_{1}(s,y(s,t,\eta)))ds\\
&-\int_{t}^{+\infty}U(t)(I-P)U^{-1}(s)f(s,y(s,t,\eta)+h_{1}(s,y(s,t,\eta)))ds\\
=&\int_{-\infty}^{+\infty}G(t,s)f(s,y(s,t,\eta)+h_{1}(s,y(s,t,\eta))ds,
\end{split}
\end{equation}
where $G(t,s)$ is defined in (\ref{for2.4}). Using the identities, we have
\begin{equation}\label{for3.6}
\begin{split}
x(t,s,x(s,\tau,\xi))=x(t,\tau,\xi), \\
y(t,s,y(s,\tau,\eta))=y(t,\tau,\eta).
\end{split}
\end{equation}
Then, if $y(t)$ is a solution of system (\ref{sys2.2}), we have
\begin{equation}\label{for3.7}
\begin{split}
h_{1}(t,y(t,\tau,\eta))=\int_{-\infty}^{+\infty}G(t,s)f(s,y(s,\tau,\eta)+h_{1}(s,y(s,\tau,\eta))ds.
\end{split}
\end{equation}
Taking
\begin{align}\label{funch}
\mathscr H_{1}(t,y(t))=y+h_{1}(t,y(t)),  t\geq0.
\end{align}
By direct differentiation (\ref{funch}), we conclude that
\begin{equation}\label{for3.8}
\begin{split}
\mathscr H_{1}'(t,y(t))=A(t)y(t)+A(t)h_{1}(t)+f(t,y(t)+h_{1}(t))=A(t)\mathscr H_{1}(t,y)+f(t,\mathscr H_{1}(t,y)).
\end{split}
\end{equation}
The above proof implies that  if $y(t)$ is a solution of (\ref{sys2.2}), then
$\mathscr H_{1}(t, y(t))$ is a solution of (\ref{sys2.1}).
Next, we show that the existence of $\mathscr L_{1}$. Set
\begin{equation}\label{for3.9}
\begin{split}
l_{1}(t,\xi)=-\int_{-\infty}^{+\infty}G(t,s)f(s,x(s,t,\xi))ds.
\end{split}
\end{equation}
Similarly to $\tilde{h}_{1}$, we can easily prove that $l_{1}\in\Omega$. From (\ref{for3.6}), if $x(t)$ is a solution of system (\ref{sys2.1}), we have
\begin{equation}\label{for3.10}
\begin{split}
l_{1}(t,x(t,\tau,\xi))=-\int_{-\infty}^{+\infty}G(t,s)f(s,x(s,\tau,\xi))ds.
\end{split}
\end{equation}
Taking
\begin{align}\label{funcl}
\mathscr L_{1}(t,x(t))=x(t)+l_{1}(t,x(t)),  t\geq0.
\end{align}
By direct differentiation (\ref{funcl}), we have
\begin{equation}\label{for3.11}
\begin{split}
\mathscr L_{1}'(t,x(t))=A(t)x(t)+f(t,x(t))+A(t)l_{1}(t,x(t))-f(t,x(t))=A(t)\mathscr L_{1}(t,x(t)).
\end{split}
\end{equation}
The above proof implies that  if $x(t)$ is a solution of (\ref{sys2.1}), then
$\mathscr L_{1}(t, x(t))$ is a solution of (\ref{sys2.2}).\\
\indent Next we prove $\mathscr H_{1}(t,\mathscr L_{1}(t,x))=x$ and $\mathscr L_{1}(t,\mathscr H_{1}(t,y))=y$, for $t\geq0$. Let $x(t)$ be any solution of system (\ref{sys2.1}). Then we know that $\mathscr L_{1}(t,x(t))\triangleq w_{1}(t)$ is the solution of system (\ref{sys2.2}) and $\mathscr H_{1}(t,\mathscr L_{1}(t,x(t)))\triangleq g_{1}(t)$ is the solution of system (\ref{sys2.1}).
It follows from (\ref{funch}) and (\ref{funcl}),
\begin{equation}\label{for3.14}
\begin{split}
&w_{1}(t)=x(t)+l_{1}(t,x(t))=x(t)-\int_{-\infty}^{+\infty}G(t,s)f(s,x(s))ds,\\
&g_{1}(t)=w_{1}(t)+h_{1}(t,w_{1}(t))=w_{1}(t)+\int_{-\infty}^{+\infty}G(t,s)f(s,g_{1}(s))ds.
\end{split}
\end{equation}
Then,
\begin{equation}\label{for3.14}
\begin{split}
&|g_{1}(t)-x(t)|\\
\leq&\int_{-\infty}^{+\infty}|G(t,s)||f(s,g_{1}(s))-f(s,x(s))|ds\\
\leq&\int_{-\infty}^{+\infty}|G(t,s)|\gamma|g_{1}(s)-x(s)|ds\\
\leq&3\beta\gamma\alpha^{-1}\|g_{1}(t)-x(t)\|.
\end{split}
\end{equation}
Therefore, $\mathscr H_{1}(t,\mathscr L_{1}(t,x))=x$.
Let $y(t)$ be any solution of system (\ref{sys2.2}). Then we know that $\mathscr H_{1}(t,y(t))\triangleq \hat{g}_{1}(t)$ is the solution of system (\ref{sys2.1}) and $\mathscr L_{1}(t,\mathscr H_{1}(t,y(t)))\triangleq \hat{w}_{1}(t)$ is the solution of system (\ref{sys2.2}).
It follows from (\ref{funch}) and (\ref{funcl}),
\begin{equation}\label{for3.16}
\begin{split}
&\hat{g}_{1}(t)=y(t)+h_{1}(t,y(t))=y(t)+\int_{-\infty}^{+\infty}G(t,s)f(s,\hat{g}_{1}(t))ds,\\
&\hat{w}_{1}(t)=\hat{g}_{1}(t)+l_{1}(t,\hat{g}_{1}(t))=\hat{g}_{1}(t)-\int_{-\infty}^{+\infty}G(t,s)f(s,\hat{g}_{1}(t))ds.
\end{split}
\end{equation}
Then, $\hat{w}_{1}(t)=y(t)$. Therefore, $\mathscr L_{1}(t,\mathscr H_{1}(t,y))=y$.\\
{\bf Step 2.} We prove if
$t<0$, Theorem \ref{thm2.6} still holds. For $t<0$, $\xi,\eta\in X$ and $v\in\Omega$, we define the following mapping:
\begin{equation}\label{for3.17}
\begin{split}
&{\mathscr F}v(t,\eta)=\tilde{v}(t,\eta)\\
=&\int_{-\infty}^{t}U(t)(I-Q)U^{-1}(s)f(s,y(s,t,\eta)+v(s,y(s,t,\eta)))ds\\
&-\int_{t}^{0}U(t)QU^{-1}(s)f(s,y(s,t,\eta)+v(s,y(s,t,\eta)))ds\\
&-\int_{0}^{+\infty}U(t)(I-P)U^{-1}(s)f(s,y(s,t,\eta)+v(s,y(s,t,\eta)))ds\\
=&\int_{-\infty}^{+\infty}G(t,s)f(s,y(s,t,\eta)+v(s,y(s,t,\eta)))ds,
\end{split}
\end{equation}
where $G(t,s)$ is defined in (\ref{for2.4}). It follows from (\ref{for2.3}), we obtain
\begin{equation}\label{for3.18}
\begin{split}
&|\tilde{v}(t,\eta)|\\
\leq& \int_{-\infty}^{t}\beta e^{-\alpha(t-s)}\mu ds
+\int_{t}^{0}\beta e^{-\alpha(t-s)}\mu ds + \int_{0}^{+\infty}\beta e^{-\alpha(s-t)}\mu ds\\
\leq&\beta\mu\alpha^{-1}+\beta\mu\alpha^{-1}+\beta\mu\alpha^{-1}\\
\leq&3\beta\mu\alpha^{-1}.
\end{split}
\end{equation}
Moreover,  by differentiating (\ref{for3.17}), we get
\begin{equation}\label{for3.19}
\begin{split}
\tilde{v}^{'}(t,\eta)=A(t)\tilde{v}(t,\eta)+f(t,\eta+v(t,\eta)).
\end{split}
\end{equation}
Hence, from (\ref{for3.17}) and (\ref{for3.18}), we know $\tilde{v}(t,\eta)$ is continuous and $\tilde{v}(t,\eta)\in \Omega$.  For any $v_{1},v_{2}\in\Omega$,
from (\ref{for2.7}), we have
\begin{equation}\label{for3.20}
\begin{split}
&|\tilde{v}_{1}(t)-\tilde{v}_{2}(t)|\\
\leq& \int_{-\infty}^{t}\beta e^{-\alpha(t-s)}\gamma|(v_{1}(s)-v_{2}(s)| ds+\int_{t}^{0}\beta e^{-\alpha(t-s)}\gamma|(v_{1}(s)-v_{2}(s)| ds\\
&+\int_{0}^{+\infty}\beta e^{-\alpha(t-s)}\gamma|(v_{1}(s)-v_{2}(s)| ds\\
\leq&3\beta\gamma\alpha^{-1}\|v_{1}-v_{2}\|.
\end{split}
\end{equation}
Note that $3\beta\gamma\alpha^{-1}<1$. Thus, $\mathscr F$ : $\Omega\to\Omega$ is a contraction map. Therefore, there exists a  unique fixed point $v\in\Omega$ such that
\begin{equation}\label{for3.21}
\begin{split}
&v(t,\eta)\\
=&\int_{-\infty}^{t}U(t)(I-Q)U^{-1}(s)f(s,y(s,t,\eta)+v(s,y(s,t,\eta)))ds\\
&-\int_{t}^{0}U(t)QU^{-1}(s)f(s,y(s,t,\eta)+v(s,y(s,t,\eta)))ds\\
&-\int_{0}^{+\infty}U(t)(I-P)U^{-1}(s)f(s,y(s,t,\eta)+v(s,y(s,t,\eta)))ds\\
\triangleq&\int_{-\infty}^{+\infty}G(t,s)f(s,y(s,t,\eta)+v(s,y(s,t,\eta))ds.
\end{split}
\end{equation}
Then, if $y(t)$ is a solution of system (\ref{sys2.2}), we have
\begin{equation}\label{for3.22}
\begin{split}
v(t,y(t,\tau,\eta))=\int_{-\infty}^{+\infty}G(t,s)f(s,y(s,\tau,\eta)+v(s,y(s,\tau,\eta))ds.
\end{split}
\end{equation}
Taking
\begin{align}\label{funcv}
\mathscr H_{2}(t,y(t))=y+v(t,y(t)), \quad t<0.
\end{align}
By direct differentiation (\ref{funcv}), we conclude that
\begin{equation}\label{for3.24}
\begin{split}
\mathscr H_{2}'(t,y(t))=A(t)y(t)+A(t)v(t)+f(t,y(t)+v(t))=A(t)\mathscr H_{2}(t,y)+f(t,\mathscr H_{2}(t,y)).
\end{split}
\end{equation}
The above proof implies that  if $y(t)$ is a solution of (\ref{sys2.2}), then
$\mathscr H_{2}(t, y(t))$ is a solution of (\ref{sys2.1}).
Next, we construct the function $\mathscr L_{2}$. Set
\begin{equation}\label{for3.25}
\begin{split}
l_{2}(t,\xi)=-\int_{-\infty}^{+\infty}G(t,s)f(s,x(s,t,\xi))ds.
\end{split}
\end{equation}
Similarly to $\tilde{v}$, we can prove that $l_{2}\in\Omega$. From (\ref{for3.6}), if $x(t)$ is a solution of system (\ref{sys2.1}), we have
\begin{equation}\label{for3.26}
\begin{split}
l_{2}(t,x(t,\tau,\xi))=-\int_{-\infty}^{+\infty}G(t,s)f(s,x(s,\tau,\xi))ds.
\end{split}
\end{equation}
Taking
\begin{align}\label{funcl2}
\mathscr  L_{2}(t,x(t))=x(t)+l_{2}(t,x(t)), \quad t<0.
\end{align}
By direct differentiation, we get
\begin{equation}\label{for3.28}
\begin{split}
\mathscr L_{2}'(t,x(t))=A(t)x(t)+f(t,x(t))+A(t)l_{2}(t,x(t))-f(t,x(t))=A(t)\mathscr  L_{2}(t,x(t)).
\end{split}
\end{equation}
The above proof implies that  if $x(t)$ is a solution of (\ref{sys2.1}), then
$\mathscr  L_{2}(t, x(t))$ is a solution of (\ref{sys2.2}).\\
\indent Next we prove $\mathscr  H_{2}(t,\mathscr  L_{2}(t,x))=x$ and $\mathscr  L_{2}(t,\mathscr H_{2}(t,y))=y$, for $t<0$. Let $x(t)$ be any solution of system (\ref{sys2.1}). Then we know that $\mathscr L_{2}(t,x(t))\triangleq w_{2}(t)$ is the solution of system (\ref{sys2.2}) and $\mathscr H_{2}(t,\mathscr L_{2}(t,x(t)))\triangleq g_{2}(t)$ is the solution of system (\ref{sys2.1}).
From (\ref{funcv}) and (\ref{funcl2}),
\begin{equation}\label{for3.29}
\begin{split}
&w_{2}(t)=x(t)+l_{2}(t,x(t))=x(t)-\int_{-\infty}^{+\infty}G(t,s)f(s,x(s))ds,\\
&g_{2}(t)=w_{2}(t)+h(t,w_{2}(t))=w_{2}(t)+\int_{-\infty}^{+\infty}G(t,s)f(s,g_{2}(s))ds.
\end{split}
\end{equation}
Then,
\begin{equation}\label{for3.30}
\begin{split}
&|g_{2}(t)-x(t)|\\
\leq&\int_{-\infty}^{+\infty}|G(t,s)||f(s,g_{2}(s))-f(s,x(s))|ds\\
\leq&\int_{-\infty}^{+\infty}|G(t,s)|\gamma|g_{2}(s)-x(s)|ds\\
\leq&3\beta\gamma\alpha^{-1}\|g_{2}(t)-x(t)\|.
\end{split}
\end{equation}
Therefore, $\mathscr H_{2}(t,\mathscr L_{2}(t,x))=x$.
Let $y(t)$ be any solution of system (\ref{sys2.2}). Then we know that $\mathscr H_{2}(t,y(t))\triangleq \hat{g}_{2}(t)$ is the solution of system (\ref{sys2.1}) and $\mathscr L_{2}(t,\mathscr H_{2}(t,y(t)))\triangleq \hat{w}_{2}(t)$ is the solution of system (\ref{sys2.2}).
From (\ref{funch}) and (\ref{funcl}),
\begin{equation}\label{for3.31}
\begin{split}
&\hat{g}_{2}(t)=y(t)+v(t,y(t))=y(t)+\int_{-\infty}^{+\infty}G(t,s)f(s,\hat{g}_{2}(t))ds,\\
&\hat{w}_{2}(t)=\hat{g}(t)+\bar{l}(t,\hat{g}_{2}(t))=\hat{g}_{2}(t)-\int_{-\infty}^{+\infty}G(t,s)f(s,\hat{g}_{2}(t))ds.
\end{split}
\end{equation}
Then, $\hat{w}_{2}(t)=y(t)$. Therefore, $\mathscr L_{2}(t,\mathscr H_{2}(t,y))=y$. We have proved that the Theorem \ref{thm2.6} holds for $t<0$. \\
\noindent{\bf Step 3.} At last, we prove that if $t\to0$, we have $H_{2}(t,y(t))\to H_{1}(0,y(0))$.  
Recall that if $t=0$, we obtain
\begin{align*}
&|h_{1}(t,y(t,\tau,\eta))-v(t,y(t,\tau,\eta))|\\
=&|\int_{-\infty}^{0}U(0)(I-Q)U^{-1}(s)\big(f(s,y(s,\tau,\eta)+h_{1}(s,y(s,\tau,\eta)))\\
&-f(s,y(s,\tau,\eta)+v(s,y(s,\tau,\eta)))\big)ds\\
&-\int_{0}^{+\infty}U(0)(I-P)U^{-1}(s)\big(f(s,y(s,\tau,\eta)+h_{1}(s,y(s,\tau,\eta)))\\
&-U^{-1}(s)f(s,y(s,\tau,\eta)+v(s,y(s,\tau,\eta)))\big)ds|\\
\leq&|\int_{-\infty}^{0}U(0)(I-Q)U^{-1}(s)\big(f(s,y(s,\tau,\eta)+h_{1}(s,y(s,\tau,\eta)))\\
&-f(s,y(s,\tau,\eta)+v(s,y(s,\tau,\eta)))\big)ds|\\
&+|\int_{0}^{+\infty}U(0)(I-P)U^{-1}(s)\big(f(s,y(s,\tau,\eta)+h_{1}(s,y(s,\tau,\eta)))\\
&-f(s,y(s,\tau,\eta)+v(s,y(s,\tau,\eta)))\big)ds|\\
\leq&\int_{-\infty}^{0}\beta e^{\alpha s}\gamma|h_{1}-v|ds+|\int_{0}^{+\infty}\beta e^{-\alpha s}\gamma|h_{1}-v|ds\\
=&2\beta\gamma\alpha^{-1}\|h_{1}-v\|.
\end{align*}
Now note that $2\beta\gamma\alpha^{-1}<1$, we have $h_{1}=v$. Thus, if $t=0$,
\begin{align*}
\mathscr H_{1}(0,y(0))=&\int_{-\infty}^{0}U(0)(I-Q)U^{-1}(s)f(s,y(s,\tau,\eta)+h_{1}(s,y(s,\tau,\eta)))ds\\
&-\int_{0}^{+\infty}U(0)(I-P)U^{-1}(s)f(s,y(s,\tau,\eta)+h_{1}(s,y(s,\tau,\eta)))ds.
\end{align*}
if $t\to0$,
\begin{align*}
\mathscr H_{2}(0,y(0))=&\int_{-\infty}^{0}U(0)(I-Q)U^{-1}(s)f(s,y(s,\tau,\eta)+v(s,y(s,\tau,\eta)))ds\\
&-\int_{0}^{+\infty}U(0)(I-P)U^{-1}(s)f(s,y(s,\tau,\eta)+v(s,y(s,\tau,\eta)))ds.
\end{align*}
Hence, if $t\to0$, we have $H_{2}(t,y(t))\to H_{1}(0,y(0))$. Therefore, the proof of Theorem \ref{thm2.6} is completed.\\
{\bf Proof  of Theorem \ref{thm2.7}.}\\
{\bf Step 1.} For $t\geq0$, we prove the equivalent function $\mathscr H_{1}$ is H\"{o}lder continuous. Suppose that constants $p>0$ and $0<q<1$. From (\ref{for2.7}), we can obtain
\begin{equation}\label{for3.32}
\begin{split}
|f(t,x_{1})-f(t,x_{2})|&=|f(t,x_{1})-f(t,x_{2})|^{1-q}|f(t,x_{1})-f(t,x_{2})|^{q}\\
&\leq2\mu^{1-q}\gamma^{q}|x_{1}-x_{2}|^{q},
\end{split}
\end{equation}
where $x_{1},x_{2}\in X$. Let $\hat{\Omega}$ denote the space of all continuous maps $\varpi$ satisfies
\begin{align*}
\hat{\Omega}:=\{\varpi\in\Omega\ |\ |\varpi(t,x_{1})-\varpi(t,x_{2})|\leq p|x_{1}-x_{2}|^{q}\},
\end{align*}
for $t\geq0$, $x_{1},x_{2}\in X$. It follows from (\ref{for3.1}), we get
\begin{equation}\label{for3.33}
\begin{split}
&{\mathscr T}\varpi(t,\eta)\\
=&\int_{-\infty}^{0}U(t)(I-Q)U^{-1}(s)f(s,y(s,t,\eta)+\varpi(s,y(s,t,\eta)))ds\\
&+\int_{0}^{t}U(t)PU^{-1}(s)f(s,y(s,t,\eta)+\varpi(s,y(s,t,\eta)))ds\\
&-\int_{t}^{+\infty}U(t)(I-P)U^{-1}(s)f(s,y(s,t,\eta)+\varpi(s,y(s,t,\eta)))ds\\
=&\int_{-\infty}^{+\infty}G(t,s)f(s,y(s,t,\eta)+\varpi(s,y(s,t,\eta))ds,
\end{split}
\end{equation}
where $\varpi \in \hat{\Omega}$. By using (\ref{for2.7}) and (\ref{for3.32})
\begin{equation}\label{for3.34}
  \begin{split}
 &| f(t,y_{1}+\omega(t,y_{1}))-f(t,y_{2}+\omega(t,y_{2}))|\\
\leq &
 \min\{ \gamma[|y_{1}-y_{2}|+|\omega(t,y_{1})-\omega(t,y_{2})|],
    2\mu^{1-q}\gamma^{q}[|y_{1}-y_{2}|+|\omega(t,y_{1})-\omega(t,y_{2})|]^{q} \} \\
\leq &
  m \min\{|y_{1}-y_{2}|+p|y_{1}-y_{2}|^{q},|y_{1}-y_{2}|+p|y_{1}-y_{2}|^{q}]^{q}\} \\
\leq &
 m\left\{
            \begin{array}{ll}
              (1+p)|y_{1}-y_{2}|^{q}, & \mathrm{if} \; |y_{1}-y_{2}|\leq 1;\;
              (\mathrm{taking}\; \mathrm{the}\; \mathrm{left}\; \mathrm{one}) \\
              (1+p)^{q}|y_{1}-y_{2}|^{q}, & \mathrm{if}\;  |y_{1}-y_{2}|> 1.\;
              (\mathrm{taking}\; \mathrm{the}\; \mathrm{right}\; \mathrm{one})
            \end{array}
          \right. \\
\leq & m(1+p)|y_{1}-y_{2}|^{q},
\end{split}
\end{equation}
where $m=\max\{\gamma,2\mu^{1-q}\gamma^{q}\}$. 
Furthermore, by using (\ref{for2.7}) and(\ref{for3.34}), we obtain
\begin{equation}\label{for3.35}
\begin{split}
&|{\mathscr T}\varpi(t,\eta_{1})-{\mathscr T}\varpi(t,\eta_{2})|\\
\leq&|\int_{-\infty}^{+\infty}G(t,s)(f(s,y(s,t,\eta_{1})+\varpi(s,y(s,t,\eta_{1})))-f(s,y(s,t,\eta_{2})+\varpi(s,y(s,t,\eta_{2}))))ds|\\
\leq&\int_{-\infty}^{+\infty}|G(t,s)|m(1+p)\Delta_{1}^{q}(t,s)|\eta_{1}-\eta_{2}|^{q}ds\\
\leq&m(1+p)|\eta_{1}-\eta_{2}|^{q}\sup_{t\in\mathbb{R}}\int_{-\infty}^{+\infty}|G(t,s)|\Delta_{1}^{q}(t,s)ds\\
\leq&p|\eta_{1}-\eta_{2}|^{q},
\end{split}
\end{equation}
for $t\geq0$, $\eta_{1},\eta_{2}\in X$. Therefore, ${\mathscr T}\varpi\in\hat{\Omega}$. Thus, the unique fixed point $h_{1}$ of ${\mathscr T}$ belongs to $\hat{\Omega}$. From ({\ref{funch}}) and ({\ref{for3.35}}), we get
\begin{equation}\label{for3.36}
\begin{split}
&|\mathscr H_{1}(t,y_{1}(t))-\mathscr H_{1}(t,y_{2}(t))|\\
\leq&|y_{1}-y_{2}|+|h_{1}(t,y_{1})-h_{1}(t,y_{2})|\\
\leq& |y_{1}-y_{2}|+p|y_{1}-y_{2}|^{q}\\
\leq& (|y_{1}-y_{2}|^{1-q}+p)|y_{1}-y_{2}|^{q}\\
\leq&\tilde{m}_{1}|y_{1}-y_{2}|^{q},
\end{split}
\end{equation}
where $\tilde{m}_{1}=1+p$, $0<|y_{1}-y_{2}|<1$. Therefore, $\mathscr H_{1}$ is H\"{o}lder continuous. Next, we prove $\mathscr L_{1}$ is also H\"{o}lder continuous. From (\ref{for3.9}), we get
\begin{equation*}
\begin{split}
l_{1}(t,\xi)=-\int_{-\infty}^{+\infty}G(t,s)f(s,x(s,t,\xi))ds.
\end{split}
\end{equation*}
Then,
\begin{equation*}
\begin{split}
l_{1}(t,\xi_{1})-l_{1}(t,\xi_{2})=-\int_{-\infty}^{+\infty}G(t,s)(f(s,x(s,t,\xi_{1}))-f(s,x(s,t,\xi_{2})))ds.
\end{split}
\end{equation*}
It follows from (\ref{for3.9}) and (\ref{for3.32}),
\begin{equation}
\begin{split}
&|l_{1}(t,\xi_{1})-l_{1}(t,\xi_{2})|\\
\leq&|\int_{-\infty}^{+\infty}G(t,s)(f(s,x(s,t,\xi_{1}))-f(s,x(s,t,\xi_{2})))ds|\\
\leq&\int_{-\infty}^{+\infty}|G(t,s)|m(1+p)\Delta_{2}^{q}(t,s)|\xi_{1}-\xi_{2}|^{q}ds\\
\leq&m(1+p)|\xi_{1}-\xi_{2}|^{q}\sup_{t\in\mathbb{R}}\int_{-\infty}^{+\infty}|G(t,s)|\Delta_{2}^{q}(t,s)ds\\
\leq&p|\xi_{1}-\xi_{2}|^{q},
\end{split}
\end{equation}
Then,
\begin{equation}\label{for3.38}
\begin{split}
&|\mathscr L_{1}(t,x_{1})-\mathscr L_{1}(t,x_{2})|\\
\leq&|x_{1}-x_{2}|+|l_{1}(t,x_{1})-l_{1}(t,x_{2})|\\
\leq& |x_{1}-x_{2}|+p|x_{1}-x_{2}|^{q}\\
\leq& (|x_{1}-x_{2}|^{1-q}+p)|x_{1}-x_{2}|^{q}\\
\leq&\tilde{m}_{2}|x_{1}-x_{2}|^{q},
\end{split}
\end{equation}
where $\tilde{m}_{2}=1+p$, $0<|x_{1}-x_{2}|<1$. Hence, for $t\geq0$, we prove $\mathscr H_{1}$ and $\mathscr L_{1}$ are H\"{o}lder continuous.\\
{\bf Step 2.} Similarly to the above proof. For $t<0$, we can easily prove the equivalent function $\mathscr H_{2}$ and $\mathscr L_{2}$ are H\"{o}lder continuous.\\
{\bf Proof of Theorem \ref{thm2.8}.}\\
Let $\Lambda$ denote the space of all continuous maps $\varrho$ satisfies
\begin{align*}
\Lambda:=\{\varrho:\mathbb{R}\to \mathbb{R}^{n}\ |\ ||\varrho|| \leq (\beta+2\kappa_{1}+\kappa_{2})c_{1}\}.
\end{align*}
Then, $(\Lambda,\|\cdot\|)$ is a Banach space. For $t\in \mathbb{R}$, $\xi,\eta\in X$ and $\varrho_{1}\in\Lambda$, we define the following mapping:
\begin{equation}\label{for3.41}
\begin{split}
&\mathscr{R} \varrho_{1}=\tilde{\varrho}_{1}\\
=&\int_{-\infty}^{t}U(t)P_{1}U^{-1}(s)f(s,y(s,t,\eta)+\varrho_{1}(s,y(s,t,\eta)))ds\\
&+\int_{0}^{t}U(t)(P_{2}+P_{3})U^{-1}(s)f(s,y(s,t,\eta)+\varrho_{1}(s,y(s,t,\eta))ds\\
&-\int_{t}^{+\infty}U(t)P_{4}U^{-1}(s)f(s,y(s,t,\eta)+\varrho_{1}(s,y(s,t,\eta)))ds\\
=&\int_{-\infty}^{+\infty}\tilde{G}(t,s)f(s,y(s,t,\eta)+\varrho_{1}(s,y(s,t,\eta)))ds,
\end{split}
\end{equation}
where $\tilde{G}(t,s)$ is defined in (\ref{foru2.7}). From (\ref{fora3.7}) and (\ref{fora3.8}), we obtain
\begin{equation}\label{fora3.17}
\begin{split}
|\int_{-\infty}^{t}U(t)P_{1}U^{-1}(s)f(s,y(s,t,\eta)+\varrho_{1}(s,y(s,t,\eta))|&\leq \int_{-\infty}^{t}\kappa_{1}\phi(s)ds<\kappa_{1}c_{1},\\
|\int_{t}^{+\infty}U(t)P_{4}U^{-1}(s)f(s,y(s,t,\eta)+\varrho_{1}(s,y(s,t,\eta))|&\leq \int_{t}^{+\infty}\kappa_{2}\phi(s)ds<\kappa_{2}c_{1}.
\end{split}
\end{equation}
From (\ref{for2.3}) and (\ref{fora3.7}), we obtain
\begin{equation}\label{fora3.18}
\begin{split}
&|\int_{0}^{t}U(t)(P_{2}+P_{3})U^{-1}(s)f(s,y(s,t,\eta)+\varrho_{1}(s,y(s,t,\eta))|\\
=&|\int_{0}^{t}U(t)PU^{-1}(s)f(s,y(s,t,\eta)+\varrho_{1}(s,y(s,t,\eta))ds\\
&-\int_{0}^{t}U(t)P_{1}U^{-1}(s)f(s,y(s,t,\eta)+\varrho_{1}(s,y(s,t,\eta))ds|\\
\leq&\int_{0}^{t}\beta e^{-\alpha(t-s)}\phi(s)ds+\int_{0}^{t}\kappa_{1}\phi(s)ds\\
\leq&(\beta+\kappa_{1})c_{1}.
\end{split}
\end{equation}
It follows from  (\ref{fora3.17}) and (\ref{fora3.18}), we get
\begin{equation}\label{for3.44}
\begin{split}
&|\mathscr{R} \varrho_{1}|\\
\leq&|\int_{-\infty}^{t}U(t)P_{1}U^{-1}(s)f(s,y(s,t,\eta)+\varrho_{1}(s,y(s,t,\eta))ds|\\
&+|\int_{0}^{t}U(t)(P_{2}+P_{3})U^{-1}(s)f(s,y(s,t,\eta)+\varrho_{1}(s,y(s,t,\eta))ds|\\
&+|\int_{t}^{+\infty}U(t)P_{4}U^{-1}(s)f(s,y(s,t,\eta)+\varrho_{1}(s,y(s,t,\eta))ds|\\
\leq&\int_{-\infty}^{0}\kappa_{1}\phi(s)ds+\int_{t}^{+\infty}\kappa_{2}\phi(s)ds+\int_{0}^{t}\beta e^{-\alpha(t-s)}\phi(s)ds+\int_{0}^{t}\kappa_{1}\phi(s)ds\\
\leq&(\beta+2\kappa_{1}+\kappa_{2})c_{1}.
\end{split}
\end{equation}
Moreover,  by differentiating (\ref{for3.41}), we get
\begin{equation}\label{for3.45}
\begin{split}
\tilde{\varrho}^{'}_{1}(t,\eta)=A(t)\tilde{\varrho}_{1}(t,\eta)+f(t,\eta+\varrho_{1}(t,\eta)).
\end{split}
\end{equation}
Hence, from (\ref{for3.44}) and (\ref{for3.45}), we know $\tilde{\varrho}_{1}(t,\eta)$ is continuous and $\tilde{\varrho}_{1}(t,\eta)\in \Lambda$.  For any $\varrho_{2},\varrho_{3}\in\Lambda$, from (\ref{fora2.8}) and (\ref{for3.41}), we have
\begin{align*}
&|\mathscr{R} \varrho_{2}-\mathscr{R} \varrho_{3}|\\
\leq&|\int_{-\infty}^{t}U(t)P_{1}U^{-1}(s)(f(s,y(s,t,\eta)+\varrho_{2}(s,y(s,t,\eta)))-f(s,y(s,t,\eta)+\varrho_{3}(s,y(s,t,\eta))))ds|\\
&+|\int_{0}^{t}U(t)(P_{2}+P_{3})U^{-1}(s)(f(s,y(s,t,\eta)+\varrho_{2}(s,y(s,t,\eta)))-f(s,y(s,t,\eta)+\varrho_{3}(s,y(s,t,\eta))))ds|\\
&+|\int_{t}^{+\infty}U(t)P_{4}U^{-1}(s)(f(s,y(s,t,\eta)+\varrho_{2}(s,y(s,t,\eta)))-f(s,y(s,t,\eta)+\varrho_{3}(s,y(s,t,\eta))))ds|\\
\leq&\int_{-\infty}^{t}\kappa_{1}\psi(s)|\varrho_{2}-\varrho_{3}|ds+\int_{0}^{t}(\beta e^{-\alpha(t-s)}+\kappa_{1})\psi(s)|\varrho_{2}-\varrho_{3}|ds+\int_{t}^{+\infty}\kappa_{2}\psi(s)|\varrho_{2}-\varrho_{3}|ds\\
\leq&(\beta+2\kappa_{1}+\kappa_{2})c_{2}\|\varrho_{2}-\varrho_{3}\|.
\end{align*}
Note that $(\beta+2\kappa_{1}+\kappa_{2})c_{2}<1$. Thus, $\mathscr{R}$ : $\Lambda\to\Lambda$ is a contraction map. Therefore, there exists a  unique fixed point $\varrho_{1}\in\Lambda$ such that
\begin{equation}\label{for3.46}
\begin{split}
&\varrho_{1}(t,\eta)\\
=&\int_{-\infty}^{t}U(t)P_{1}U^{-1}(s)f(s,y(s,t,\eta)+\varrho_{1}(s,y(s,t,\eta)))ds\\
&+\int_{0}^{t}U(t)(P_{2}+P_{3})U^{-1}(s)f(s,y(s,t,\eta)+\varrho_{1}(s,y(s,t,\eta)))ds\\
&-\int_{t}^{+\infty}U(t)P_{4}U^{-1}(s)f(s,y(s,t,\eta)+\varrho_{1}(s,y(s,t,\eta)))ds\\
=&\int_{-\infty}^{+\infty}\tilde{G}(t,s)f(s,y(s,t,\eta)+\varrho_{1}(s,y(s,t,\eta)))ds,
\end{split}
\end{equation}
where $\tilde{G}(t,s)$ is defined in (\ref{foru2.7}). By using  identities (\ref{for3.6}), if $y(t)$ is a solution of system (\ref{sys2.2}), we have
\begin{equation}\label{for3.47}
\begin{split}
&\varrho_{1}(t,y(t,\tau,\eta))\\
=&\int_{-\infty}^{t}U(t)P_{1}U^{-1}(s)f(s,y(s,\tau,\eta)+\varrho_{1}(s,y(s,\tau,\eta)))ds\\
&+\int_{0}^{t}U(t)(P_{2}+P_{3})U^{-1}(s)f(s,y(s,\tau,\eta)+\varrho_{1}(s,y(s,\tau,\eta)))ds\\
&-\int_{t}^{+\infty}U(t)P_{4}U^{-1}(s)f(s,y(s,\tau,\eta)+\varrho_{1}(s,y(s,\tau,\eta)))ds\\
=&\int_{-\infty}^{+\infty}\tilde{G}(t,s)f(s,y(s,\tau,\eta)+\varrho_{1}(s,y(s,\tau,\eta)))ds,
\end{split}
\end{equation}
where $\tilde{G}(t,s)$ is defined in (\ref{foru2.7}). Taking
\begin{align}\label{funch3}
\mathscr H_{3}(t,y(t))=y(t)+\varrho_{1}(t,y(t)),
\end{align}
By direct differentiation (\ref{funch3}), we conclude that
\begin{equation}\label{for3.48}
\begin{split}
&\mathscr H_{3}'(t,y(t))\\
=&A(t)y(t)+A(t)\varrho_{1}(t)+f(t,y(t)+\varrho_{1}(t))\\
=&A(t)\mathscr H_{3}(t,y(t))+f(t,\mathscr H_{3}(t,y(t))).
\end{split}
\end{equation}
The above proof implies that  if $y(t)$ is a solution of (\ref{sys2.2}), then
$\mathscr H_{3}(t, y(t))$ is a solution of (\ref{sys2.1}).
Next, we show that the existence of $\mathscr L_{3}$. Set
\begin{equation}\label{for3.49}
\begin{split}
&\rho_{1}(t,\xi)\\
=&-\int_{-\infty}^{t}U(t)P_{1}U^{-1}(s)f(s,x(s,t,\xi))ds\\
&-\int_{0}^{t}U(t)(P_{2}+P_{3})U^{-1}(s)f(s,x(s,t,\xi))ds\\
&+\int_{t}^{+\infty}U(t)P_{4}U^{-1}(s)f(s,x(s,t,\xi))ds.
\end{split}
\end{equation}
Similarly to $\tilde{\varrho}_{1}$, we can easily prove that $\rho_{1}\in\Lambda$. From (\ref{for3.6}), if $x(t)$ is a solution of system (\ref{sys2.1}), we have
\begin{equation}\label{for3.51}
\begin{split}
&\rho_{1}(t,x(t,\tau,\xi))\\
=&-\int_{-\infty}^{t}U(t)P_{1}U^{-1}(s)f(s,x(s,\tau,\xi))ds\\
&-\int_{0}^{t}U(t)(P_{2}+P_{3})U^{-1}(s)f(s,x(s,\tau,\xi))ds\\
&+\int_{t}^{+\infty}U(t)P_{4}U^{-1}(s)f(s,x(s,\tau,\xi))ds.
\end{split}
\end{equation}
Taking
\begin{align}\label{funcl3}
\mathscr L_{3}(t,x(t))=x(t)+\rho_{1}(t,x(t)).
\end{align}
By direct differentiation (\ref{funcl}), we have
\begin{equation}\label{for3.53}
\begin{split}
&\mathscr L_{3}'(t,x(t))\\
=&A(t)x(t)+f(t,x(t))+A(t)\rho_{1}(t,x(t))-f(t,x(t))\\
=&A(t)\mathscr L_{3}(t,x(t)).
\end{split}
\end{equation}
The above proof implies that  if $x(t)$ is a solution of (\ref{sys2.1}), then
$\mathscr L_{3}(t, x(t))$ is a solution of (\ref{sys2.2}).\\
\indent Next we prove $\mathscr H_{3}(t,\mathscr L_{3}(t,x))=x$ and $\mathscr L_{3}(t,\mathscr H_{3}(t,y))=y$, for $t\in\mathbb{R}$. Let $x(t)$ be any solution of system (\ref{sys2.1}). Then we know that $\mathscr L_{3}(t,x(t))\triangleq \omega_{1}(t)$ is the solution of system (\ref{sys2.2}) and $\mathscr H_{3}(t,\mathscr L_{3}(t,x(t)))\triangleq \omega_{2}(t)$ is the solution of system (\ref{sys2.1}).
\begin{equation}\label{for3.54}
\begin{split}
&\omega_{1}(t)=x(t)+\rho_{1}(t,x(t))\\
=&x(t)-\int_{-\infty}^{t}U(t)P_{1}U^{-1}(s)f(s,x(s))ds-\int_{0}^{t}U(t)(P_{2}+P_{3})U^{-1}(s)f(s,x(s))ds\\
&+\int_{t}^{+\infty}U(t)P_{4}U^{-1}(s)f(s,x(s))ds,\\
&\omega_{2}(t)=\omega_{1}(t)+\varrho_{1}(t,\omega_{1}(t))\\
=&\omega_{1}(t)+\int_{-\infty}^{t}U(t)P_{1}U^{-1}(s)f(s,\omega_{2}(s))ds+\int_{0}^{t}U(t)(P_{2}+P_{3})U^{-1}(s)f(s,\omega_{2}(s))ds\\
&-\int_{t}^{+\infty}U(t)P_{4}U^{-1}(s)f(s,\omega_{2}(s))ds.
\end{split}
\end{equation}
Then,
\begin{equation}\label{for3.14}
\begin{split}
&|\omega_{2}(t)-x(t)|\\
\leq&\int_{-\infty}^{t}\kappa_{1}\psi(s)|\omega_{2}(s)-x(s)|ds+\int_{0}^{t}(\beta e^{-\alpha(t-s)}+\kappa_{1})\psi(s)|\omega_{2}(s)-x(s)|ds\\
&+\int_{t}^{+\infty}\kappa_{2}\psi(s)|\omega_{2}(s)-x(s)|ds\\
\leq&(\beta+2\kappa_{1}+\kappa_{2})c_{2}\|\omega_{2}(t)-x(t)\|.
\end{split}
\end{equation}
Therefore, $\omega_{2}(t)=x(t)$. Let $y(t)$ be any solution of system (\ref{sys2.2}). Then we know that $\mathscr H_{3}(t,y(t))\triangleq \hat{\omega}_{1}(t)$ is the solution of system (\ref{sys2.1}) and $\mathscr L_{3}(t,\mathscr H_{3}(t,y(t)))\triangleq \hat{\omega}_{2}(t)$ is the solution of system (\ref{sys2.2}).
\begin{equation}\label{for3.16}
\begin{split}
&\hat{\omega}_{1}(t)=y(t)+\varrho_{1}(t,y(t))\\
=&y(t)+\int_{-\infty}^{t}U(t)P_{1}U^{-1}(s)f(s,\hat{\omega}_{1}(s))ds+\int_{0}^{t}U(t)(P_{2}+P_{3})U^{-1}(s)f(s,\hat{\omega}_{1}(s))ds\\
&-\int_{t}^{+\infty}U(t)P_{4}U^{-1}(s)f(s,\hat{\omega}_{1}(s))ds,\\
&\hat{\omega}_{2}(t)=\hat{\omega}_{1}(t)+\rho_{1}(t,\hat{\omega}_{1}(t))\\
=&\hat{\omega}_{1}(t)-\int_{-\infty}^{t}U(t)P_{1}U^{-1}(s)f(s,\hat{\omega}_{1}(s))ds-\int_{0}^{t}U(t)(P_{2}+P_{3})U^{-1}(s)f(s,\hat{\omega}_{1}(s))ds\\
&+\int_{t}^{+\infty}U(t)P_{4}U^{-1}(s)f(s,\hat{\omega}_{1}(s))ds.
\end{split}
\end{equation}
Then, $\hat{\omega}_{2}(t)=y(t)$. Therefore, $\mathscr L_{3}(t,\mathscr H_{3}(t,y))=y$. Hence, the proof of Theorem \ref{thm2.8} is completed.\\
{\bf Proof  of Theorem \ref{thmh}.}\\
{\bf Step 1.}We prove the equivalent function $\mathscr H_{1}$ is H\"{o}lder continuous. Suppose that constants $c\geq1$, $p>0$ and $0<q<1$. From (\ref{fora2.8}), we can obtain
\begin{equation}\label{for3.57}
\begin{split}
|f(t,x_{1})-f(t,x_{2})|&=|f(t,x_{1})-f(t,x_{2})|^{1-q}|f(t,x_{1})-f(t,x_{2})|^{q}\\
&\leq2\phi^{1-q}(t)\psi^{q}(t)|x_{1}-x_{2}|^{q},
\end{split}
\end{equation}
where $x_{1},x_{2}\in X$. Let $\hat{\Lambda}$ denote the space of all continuous maps $\chi$ satisfies
\begin{align*}
\hat{\Lambda}:=\{\chi\in\Lambda\ |\ |\chi(t,x_{1})-\chi(t,x_{2})|\leq p|x_{1}-x_{2}|^{q}\},
\end{align*}
for $t\in \mathbb{R}$, $x_{1},x_{2}\in X$. It follows from (\ref{for3.41}), we get
\begin{equation}\label{for3.58}
\begin{split}
&{\mathscr R}\chi(t,\eta)\\
=&\int_{-\infty}^{t}U(t)P_{1}U^{-1}(s)f(s,y(s,t,\eta)+\chi(s,y(s,t,\eta)))ds\\
&+\int_{0}^{t}U(t)(P_{2}+P_{3})U^{-1}(s)f(s,y(s,t,\eta)+\chi(s,y(s,t,\eta))ds\\
&-\int_{t}^{+\infty}U(t)P_{4}U^{-1}(s)f(s,y(s,t,\eta)+\chi(s,y(s,t,\eta)))ds\\
=&\int_{-\infty}^{+\infty}\tilde{G}(t,s)f(s,y(s,t,\eta)+\chi(s,y(s,t,\eta))ds,
\end{split}
\end{equation}
where $\tilde{G}(t,s)$ is defined in (\ref{foru2.7}). By using (\ref{fora2.8}) and (\ref{for3.57})
\begin{equation}\label{for3.59}
  \begin{split}
 &| f(t,y_{1}+\chi(t,y_{1}))-f(t,y_{2}+\chi(t,y_{2}))|\\
\leq &
 \min\{ \psi(t)[|y_{1}-y_{2}|+|\chi(t,y_{1})-\chi(t,y_{2})|],
    \phi^{1-q}(t)\psi^{q}(t)[|y_{1}-y_{2}|+|\chi(t,y_{1})-\chi(t,y_{2})|]^{q} \} \\
\leq &
  M(t)\min\{|y_{1}-y_{2}|+p|y_{1}-y_{2}|^{q},|y_{1}-y_{2}|+p|y_{1}-y_{2}|^{q}]^{q}\} \\
\leq &
  M(t)\left\{
            \begin{array}{ll}
              (1+p)|y_{1}-y_{2}|^{q}, & \mathrm{if} \; |y_{1}-y_{2}|\leq 1;\;
              (\mathrm{taking}\; \mathrm{the}\; \mathrm{left}\; \mathrm{one}) \\
              (1+p)^{q}|y_{1}-y_{2}|^{q}, & \mathrm{if}\;  |y_{1}-y_{2}|> 1.\;
              (\mathrm{taking}\; \mathrm{the}\; \mathrm{right}\; \mathrm{one})
            \end{array}
          \right. \\
\leq & M(t)(1+p)|y_{1}-y_{2}|^{q}.
\end{split}
\end{equation}
Furthermore, by using (\ref{fora2.8}) and (\ref{for3.59}), we obtain
\begin{equation}\label{for3.60}
\begin{split}
&|{\mathscr R}\chi(t,\eta_{1})-{\mathscr R}\chi(t,\eta_{2})|\\
\leq&\int_{-\infty}^{+\infty}\tilde{G}(t,s)(f(s,y(s,t,\eta_{1})+\chi(s,y(s,t,\eta_{1})))-f(s,y(s,t,\eta_{2})+\chi(s,y(s,t,\eta_{2}))))ds\\
\leq&\int_{-\infty}^{+\infty}|\tilde{G}(t,s)|M(s)(1+p)\Delta_{1}^{q}(t,s)|\eta_{1}-\eta_{2}|^{q}ds\\
\leq&(1+p)|\eta_{1}-\eta_{2}|^{q}\sup_{t\in\mathbb{R}}\int_{-\infty}^{+\infty}|\tilde{G}(t,s)|M(s)\Delta_{1}^{q}(t,s)ds\\
\leq&p|\eta_{1}-\eta_{2}|^{q},
\end{split}
\end{equation}
for $t\in\mathbb{R}$, $\eta_{1},\eta_{2}\in X$, where $\tilde{G}(t,s)$ is defined in (\ref{foru2.7}). Therefore, ${\mathscr R}\chi\in\hat{\Lambda}$. Thus, the unique fixed point $\varrho_{1}$ of ${\mathscr R}$ belongs to $\hat{\Lambda}$. From ({\ref{funch}}) and ({\ref{for3.35}}), we get
\begin{equation}\label{for3.61}
\begin{split}
&|\mathscr H_{3}(t,y_{1}(t))-\mathscr H_{3}(t,y_{2}(t))|\\
\leq&|y_{1}-y_{2}|+|\varrho_{1}(t,y_{1})-\varrho_{1}(t,y_{2})|\\
\leq& |y_{1}-y_{2}|+p|y_{1}-y_{2}|^{q}\\
\leq& (|y_{1}-y_{2}|^{1-q}+p)|y_{1}-y_{2}|^{q}\\
\leq&\tilde{c}|y_{1}-y_{2}|^{q},
\end{split}
\end{equation}
where $\tilde{c}=1+p$, $0<|y_{1}-y_{2}|<1$. Therefore, $\mathscr H_{3}$ is H\"{o}lder continuous. Next, we prove $\mathscr L_{3}$ is also H\"{o}lder continuous. From (\ref{for3.49}), we get
\begin{equation*}
\begin{split}
&\rho_{1}(t,\xi)\\
=&-\int_{-\infty}^{t}U(t)P_{1}U^{-1}(s)f(s,x(s,t,\xi))ds\\
&-\int_{0}^{t}U(t)(P_{2}+P_{3})U^{-1}(s)f(s,x(s,t,\xi))ds\\
&+\int_{t}^{+\infty}U(t)P_{4}U^{-1}(s)f(s,x(s,t,\xi))ds\\
=&-\int_{-\infty}^{+\infty}\tilde{G}(t,s)f(s,x(s,t,\xi))ds.
\end{split}
\end{equation*}
where $\tilde{G}(t,s)$ is defined in (\ref{foru2.7}). Then,
\begin{equation*}
\begin{split}
&\rho_{1}(t,\xi_{1})-\rho_{1}(t,\xi_{2})\\
=&-\int_{-\infty}^{+\infty}\tilde{G}(t,s)(f(s,x(s,t,\xi_{1}))-f(s,x(s,t,\xi_{2})))ds.
\end{split}
\end{equation*}
It follows from (\ref{for3.49}) and (\ref{for3.57}),
\begin{equation}\label{for3.62}
\begin{split}
&|\rho_{1}(t,\xi_{1})-\rho_{1}(t,\xi_{2})|\\
\leq&\int_{-\infty}^{+\infty}\tilde{G}(t,s)(f(s,x(s,t,\xi_{1}))-f(s,x(s,t,\xi_{2})))ds\\
\leq&\int_{-\infty}^{+\infty}|\tilde{G}(t,s)|(1+p)M(s)\Delta_{2}^{q}(t,s)|\xi_{1}-\xi_{2}|^{q}ds\\
\leq&(1+p)|\xi_{1}-\xi_{2}|^{q}\sup_{t\in\mathbb{R}}\int_{-\infty}^{+\infty}|\tilde{G}(t,s)|M(s)\Delta_{2}^{q}(t,s)ds\\
\leq&p|\xi_{1}-\xi_{2}|^{q},
\end{split}
\end{equation}
Then,
\begin{equation}\label{for3.63}
\begin{split}
&|\mathscr L_{3}(t,x_{1})-\mathscr L_{3}(t,x_{2})|\\
\leq&|x_{1}-x_{2}|+|\rho_{1}(t,x_{1})-\rho_{1}(t,x_{2})|\\
\leq& |x_{1}-x_{2}|+p|x_{1}-x_{2}|^{q}\\
\leq& (|x_{1}-x_{2}|^{1-q}+p)|x_{1}-x_{2}|^{q}\\
\leq&\tilde{c}|x_{1}-x_{2}|^{q},
\end{split}
\end{equation}
where $\tilde{c}=1+p$, $0<|x_{1}-x_{2}|<1$. Hence, we prove $\mathscr H_{3}$ and $\mathscr L_{3}$ are H\"{o}lder continuous.\\
{\bf Proof of Theorem {\ref{thm2.10}}.} Firstly, we prove some lemmas, consider periodic system
\begin{equation}\label{equa4.1}
\begin{split}
x'=\sigma(t,x),
\end{split}
\end{equation}
where $\sigma(t+T,x)=\sigma(t,x)$. Systems (\ref{equa4.1}) satisfy the existence and uniqueness of the solution. Suppose that $\bar{X}(t,(t_{0}, x_{0}))$ is the solution of system (\ref{equa4.1}) satisfying
$$\bar{X}(t_{0},(t_{0}, \bar{x}_{0}))=\bar{x}_{0}.$$
\begin{lemma}\label{lem4.1}
For any $t,s\in \mathbb{R}$, $x\in\mathbb{R}^{n}$, we have
\begin{align*}
\bar{X}(t+T,(s+T, x))&=\bar{X}(t,(s, x)).
\end{align*}
\end{lemma}
{\bf Proof}. From variation formula, we have
\begin{align}
\bar{X}(t,(s, x))&=x+\int_{s}^{t}\sigma(\tau,\bar{X}(\tau,(s, x)))d\tau.
\end{align}
Then,
\begin{align}
&\bar{X}(t+T,(s+T,x))\nonumber\\
&=x+\int_{s+T}^{t+T}\sigma(\tau,\bar{X}(\tau,(s+T, x)))d\tau\nonumber\\
&\overset{\tau=\tau_{1}+T}{=}x+\int_{s}^{t}\sigma(\tau_{1},\bar{X}(\tau_{1}+T,(s+T, x)))d\tau_{1}.\label{for4.3}
\end{align}
Denote $F_{1}(t)=\bar{X}(t+T,(s+T, x))$. By (\ref{for4.3}), we know $F_{1}(t)$ is the solution of system (\ref{equa4.1}). Since $F_{1}(s)=x$,
$\bar{X}(s,(s,x))=x$, according to the existence and uniqueness of the solution, we get $\bar{X}(t+T,(s+T, x))=\bar{X}(t,(s, x))$.

\begin{lemma}\label{lem4.2}
Suppose that  periodic system $x'=A(t)x$ $(A(t+T)=A(t))$ have an exponential trichotomy ($U(t)$ is the fundamental matrix of system (\ref{sys2.2}) satisfying (\ref{for2.3})).  For any $t,s\in \mathbb{R}$, we have
\begin{align*}
G(t+T,s+T)=G(t,s),\quad \tilde{G}(t+T,s+T)=\tilde{G}(t,s).
\end{align*}
\end{lemma}
{\bf Proof}. $U(t)$ is the fundamental matrix  of linear system $x'=A(t)x$. It's easy to get $U(t+T)$ is also the fundamental matrix of linear system $x'=A(t)x$. Thus, there exists an invertible matrix $D$ such that $U(t+T) =U(t)D$. Taking $B=\frac{1}{T}\ln D$ and $M(t)=U(t)e^{-Bt}$. Then,
\begin{equation}\label{for4.4}
\begin{split}
M(t+T)&=U(t+T)e^{-B(t+T)}=U(t)DD^{-1}e^{-Bt}=M(t).\\
M^{-1}(t+T)&=e^{B(t+T)}U^{-1}(t+T)=e^{Bt}DD^{-1}U^{-1}(t)=M^{-1}(t).
\end{split}
\end{equation}
From (\ref{for4.4}), we get
\begin{align*}
U(t+T)P_{1}U^{-1}(s+T)&=M(t+T)e^{B(t+T)}P_{1}e^{-B(s+T)}M^{-1}(s+T)\\
&=M(t)e^{Bt}P_{1}e^{-Bs}M^{-1}(s)\\
&=U(t)P_{1}U^{-1}(s).
\end{align*}
Similar to the above proof, we can easily get $U(t+T)(P_{2}+P_{3})U^{-1}(s+T)=U(t)(P_{2}+P_{3})U^{-1}(s)$ and $U(t+T)P_{4}U^{-1}(s+T)=U(t)P_{4}U^{-1}(s)$. Thus, $\tilde{G}(t+T,s+T)=\tilde{G}(t,s)$. Similar to $\tilde{G}$, we get $G(t+T,s+T)=G(t,s)$.

Now we start to prove the periodicity or asymptotically periodic of $\mathscr H(t,x)$ and $\mathscr L(t,y)$. From   (\ref{funch3}), (\ref{funcl3}), Lemma \ref{lem4.1} and Lemma \ref{lem4.2}, we obtain,
\begin{align*}
&\mathscr H_{3}(t+T,y)\\
=&y+\int_{-\infty}^{t+T}U(t+T)P_{1}U^{-1}(s)f(s,y(s,t+T,\eta)+\varrho_{1}(s))ds\\
&+\int_{0}^{t+T}U(t+T)(P_{2}+P_{3})U^{-1}(s)f(s,y(s,t+T,\eta)+\varrho_{1}(s))ds\\
&-\int_{t+T}^{+\infty}U(t+T)P_{4}U^{-1}(s)f(s,y(s,t+T,\eta)+\varrho_{1}(s))ds\\
\overset{s=s_{1}+T}{=}&y+\int_{-\infty}^{t}U(t+T)P_{1}U^{-1}(s_{1}+T)f(s_{1}+T,y(s_{1}+T,t+T,\eta)+\varrho_{1}(s_{1}+T))ds_{1}\\
&+\int_{-T}^{t}U(t+T)(P_{2}+P_{3})U^{-1}(s_{1}+T)f(s_{1}+T,y(s_{1}+T,t+T,\eta)+\varrho_{1}(s_{1}+T))ds_{1}\\
&-\int_{t}^{+\infty}U(t+T)P_{4}U^{-1}(s_{1}+T)f(s_{1}+T,y(s_{1}+T,t+T,\eta)+\varrho_{1}(s_{1}+T))ds_{1}\\
=&y+\int_{-\infty}^{t}U(t)P_{1}U^{-1}(s_{1})f(s_{1},y(s_{1},t,\eta)+\varrho_{1}(s_{1}))ds_{1}\\
&+\int_{-T}^{t}U(t)(P_{2}+P_{3})U^{-1}(s_{1})f(s_{1},y(s_{1},t,\eta)+\varrho_{1}(s_{1}))ds_{1}\\
&-\int_{t}^{+\infty}U(t)P_{4}U^{-1}(s_{1})f(s_{1},y(s_{1},t,\eta)+\varrho_{1}(s_{1}))ds_{1}\\
\neq&\mathscr H_{3}(t,y).
\end{align*}
In addition,
\begin{align*}
&\mathscr L_{3}(t+T,x)\\
=&x-\int_{-\infty}^{t+T}U(t+T)P_{1}U^{-1}(s)f(s,x(s,t+T,\xi))ds\\
&-\int_{0}^{t+T}U(t+T)(P_{2}+P_{3})U^{-1}(s)f(s,x(s,t+T,\xi))ds\\
&+\int_{t+T}^{+\infty}U(t+T)P_{4}U^{-1}(s)f(s,x(s,t+T,\xi))ds\\
\overset{s=s_{1}+T}{=}&x-\int_{-\infty}^{t}U(t+T)P_{1}U^{-1}(s_{1}+T)f(s_{1}+T,x(s_{1}+T,t+T,\xi))ds_{1}\\
&-\int_{-T}^{t}U(t+T)(P_{2}+P_{3})U^{-1}(s_{1}+T)f(s_{1}+T,x(s_{1}+T,t+T,\xi))ds_{1}\\
&+\int_{t}^{+\infty}U(t+T)P_{4}U^{-1}(s_{1}+T)f(s_{1}+T,x(s_{1}+T,t+T,\xi))ds_{1}\\
=&x-\int_{-\infty}^{t}U(t)P_{1}U^{-1}(s_{1})f(s_{1},x(s_{1},t,\xi))ds_{1}\\
&-\int_{-T}^{t}U(t)(P_{2}+P_{3})U^{-1}(s_{1})f(s_{1},x(s_{1},t,\xi))ds_{1}\\
&+\int_{t}^{+\infty}U(t)P_{4}U^{-1}(s_{1})f(s_{1},x(s_{1},t,\xi))ds_{1}\\
\neq&\mathscr L_{3}(t,x).
\end{align*}
Thus, the equivalent function $\mathscr H_{3}$ and its inverse $\mathscr L_{3}$ in theorem \ref{thm2.8} do not have periodicity.

Next we  prove the equivalent function $\mathscr H_{3}$ and its inverse $\mathscr L_{3}$ in theorem \ref{thm2.8} do not have asymptotical periodicity.
Since
\begin{align*}
&\mathscr H_{3}(t+T,y)\\
=&y+\int_{-\infty}^{t}U(t)P_{1}U^{-1}(s)f(s,y(s,t,\eta)+\varrho_{1}(s))ds\\
&+\int_{-T}^{t}U(t)(P_{2}+P_{3})U^{-1}(s)f(s,y(s,t,\eta)+\varrho_{1}(s))ds\\
&-\int_{t}^{+\infty}U(t)P_{4}U^{-1}(s)f(s_{1},y(s,t,\eta)+\varrho_{1}(s))ds\\
=&\mathscr H_{3}(t,y)+\int_{-T}^{0}U(t)(P_{2}+P_{3})U^{-1}(s)f(s,y(s,t,\eta)+\varrho_{1}(s))ds
\end{align*}
and
\begin{align*}
\lim_{t\to\infty}\int_{-T}^{0}U(t)(P_{2}+P_{3})U^{-1}(s)f(s,y(s,t,\eta)+\varrho_{1}(s))ds\neq0.
\end{align*}
Thus, $\mathscr H_{3}$ do not have asymptotical periodicity. The proof of $\mathscr L_{3}$ is similar to $\mathscr H_{3}$, we omit.

Next, we prove equivalent function $\mathscr H_{1}$, $\mathscr H_{2}$, $\mathscr L_{1}$ and $\mathscr L_{2}$ in theorem \ref{thm2.6} do not have periodicity or asymptotical periodicity.
\begin{align*}
&\mathscr H_{1}(t+T,y)\\
=&y+\int_{-\infty}^{0}U(t+T)(I-Q)U^{-1}(s)f(s,y(s,t+T,\eta)+h_{1}(s))ds\\
&+\int_{0}^{t+T}U(t+T)PU^{-1}(s)f(s,y(s,t+T,\eta)+h_{1}(s))ds\\
&-\int_{t+T}^{+\infty}U(t+T)(I-P)U^{-1}(s)f(s,y(s,t+T,\eta)+h_{1}(s))ds\\
\overset{s=s_{1}+T}{=}&y+\int_{-\infty}^{-T}U(t+T)(I-Q)U^{-1}(s_{1}+T)f(s_{1}+T,y(s_{1}+T,t+T,\eta)+h_{1}(s_{1}+T))ds_{1}\\
&+\int_{-T}^{t}U(t+T)PU^{-1}(s_{1}+T)f(s_{1}+T,y(s_{1}+T,t+T,\eta)+h_{1}(s_{1}+T))ds_{1}\\
&-\int_{t}^{+\infty}U(t+T)(I-P)U^{-1}(s_{1}+T)f(s_{1}+T,y(s_{1}+T,t+T,\eta)+h_{1}(s_{1}+T))ds_{1}\\
=&y+\int_{-\infty}^{-T}U(t)(I-Q)U^{-1}(s_{1})f(s_{1},y(s_{1},t,\eta)+h_{1}(s_{1}))ds_{1}\\
&+\int_{-T}^{t}U(t)PU^{-1}(s_{1})f(s_{1},y(s_{1},t,\eta)+h_{1}(s_{1}))ds_{1}\\
&-\int_{t}^{+\infty}U(t)(I-P)U^{-1}(s_{1})f(s_{1}+T,y(s_{1},t,\eta)+h_{1}(s_{1}))ds_{1}\\
\neq& \mathscr H_{1}(t,y)
\end{align*}
and
\begin{align*}
&\mathscr H_{1}(t+T,y)\\
=&y+\int_{-\infty}^{-T}U(t)(I-Q)U^{-1}(s)f(s,y(s,t,\eta)+h_{1}(s))ds\\
&+\int_{-T}^{t}U(t)PU^{-1}(s)f(s,y(s,t,\eta)+h_{1}(s))ds\\
&-\int_{t}^{+\infty}U(t)(I-P)U^{-1}(s)f(s+T,y(s,t,\eta)+h_{1}(s))ds\\
=&\mathscr H_{1}(t,y)+\int_{-T}^{0}U(t)PU^{-1}(s)f(s,y(s,t,\eta)+h_{1}(s))ds\\
&+\int_{0}^{-T}U(t)(I-Q)U^{-1}(s)f(s,y(s,t,\eta)+h_{1}(s))ds
\end{align*}
and
\begin{align*}
\lim_{t\to\infty}\big(\int_{-T}^{0}U(t)(P_{2}+P_{3})U^{-1}(s)f(s,y(s,t,\eta)+\varrho_{1}(s))ds\\
+\int_{0}^{-T}U(t)(I-Q)U^{-1}(s)f(s,y(s,t,\eta)+h_{1}(s))ds\big)\neq0.
\end{align*}
Hence, $\mathscr H_{1}$  do not have periodicity or asymptotically periodicity. The almost periodicity is similar to the above proof, we omit. The proof of $\mathscr H_{2}$, $\mathscr L_{1}$ and $\mathscr L_{2}$ are similar to $\mathscr H_{1}$, we omit.

\section{Some example}
\begin{exam}
\rm{Consider the scalar equations}
\begin{equation}\label{exam5.1}
x'=(-\frac{e^{t}-e^{-t}}{e^{t}+e^{-t}})x
\end{equation}
and
\begin{equation}\label{exam5.2}
x'=(-\frac{e^{t}-e^{-t}}{e^{t}+e^{-t}})x+\delta\sin t\sin x(t),
\end{equation}
where $\delta$ is a  positive constant.
\end{exam}
From example \ref{exam2.5}, we know that  equation (\ref{equ5.1}) admits an exponential trichotomy with $\beta=2, \alpha=1$. Let, $f(t,x)=\delta\sin t\sin x(t)$, then
\begin{align*}
&|f(t,x)|\leq\delta,\\
&|f(t,x_{1})-f(t,x_{2})|\leq\delta|x_{1}-x_{2}|.
\end{align*}
Hence, equation (\ref{exam5.2}) satisfies the condition of Theorem \ref{thm2.6} if $\delta<\frac{1}{6}$. Therefore, equation (\ref{exam5.1}) is topologically conjugated to equation (\ref{exam5.2}).

\begin{exam}
\rm{Consider the scalar equations}
\begin{equation}\label{equ5.1}
x'=(-\frac{e^{t}-e^{-t}}{e^{t}+e^{-t}})x
\end{equation}
and
\begin{equation}\label{equ5.2}
x'=(-\frac{e^{t}-e^{-t}}{e^{t}+e^{-t}})x+\frac{\epsilon }{(1+t^{2})}\sin x,
\end{equation}
where $\epsilon$ is a sufficiently small positive constant.
\end{exam}
From example \ref{exam2.5}, we know that  equation (\ref{equ5.1}) admits an exponential trichotomy with $\beta=2, \alpha=1$. Then
$U(t)=\frac{2}{e^{t}+e^{-t}}$ is a solution of (\ref{equ5.1}). Taking $P=1$, we get $|U(t)PU^{-1}(s)|<1$ $(t\geq s)$. Apart from this, we have
$|\displaystyle \int_{-\infty}^{+\infty}\frac{\epsilon}{1+t^{2}}\sin x dt|<\displaystyle\int_{-\infty}^{+\infty}\frac{\epsilon}{1+t^{2}}dt<\epsilon\pi$. Furthmore,
\begin{align*}
|\frac{\epsilon}{1+t^{2}}\sin x_{1} dt-\frac{\epsilon}{1+t^{2}}\sin x_{2} dt|<\frac{\epsilon}{1+t^{2}}|x_{1}-x_{2}|.
\end{align*}
Hence, equation (\ref{equ5.2}) satisfies the condition of Theorem \ref{thm2.8} if $\epsilon<\frac{1}{4\pi}$. Therefore, equation (\ref{equ5.1}) is topologically conjugated to equation (\ref{equ5.2}).
\section*{Data Availability Statement}
\hskip\parindent
No data was used for the research in this article. It is pure mathematics.

\section*{Conflict of Interest}
\hskip\parindent
The authors declare that they have no conflict of interest.

\section*{Contributions}
\hskip\parindent
 We declare that all the authors have same contributions to this paper.

\section*{Ethical Approval}
Not applicable.

\section*{Funding}
This paper was jointly supported from the Natural Science Foundation of Zhejiang Province (No. LZ23A010001), National Natural
Science Foundation of China under Grant (No. 11671176, 11931016), Grant Fondecyt (No. 1170466), Fondecyt 038-2021-Per\'u.

\end{document}